\documentclass[10pt, reqno]{amsart}

\usepackage{amsmath,amscd, graphics}
\usepackage{amsfonts,amssymb}
\usepackage{enumitem}

\theoremstyle{plain}
\newtheorem*{theorem*}{Theorem}
\newtheorem*{lemma*} {Lemma}
\newtheorem*{corollary*} {Corollary}
\newtheorem*{proposition*} {Proposition}
\newtheorem{theorem}{Theorem}[section]
\newtheorem{lemma}[theorem]{Lemma}
\newtheorem{corollary}[theorem]{Corollary}
\newtheorem{proposition}[theorem]{Proposition}

\theoremstyle{definition}
\newtheorem{defn}[theorem]{Definition}

\def \bn{\begin{enumerate}}
\def \en{\end{enumerate}}
\def \bdm{\begin{displaymath}}
\def \edm{\end{displaymath}}

\def \bp{\begin{proof}}
\def\ep{\end{proof}}
\def\be{\begin{equation}}
\def\ee{\end{equation}}

\begin{document}
\title{Vector fields on projective Stiefel manifolds\\ and the Browder-Dupont invariant}
\title[Vector fields on projective Stiefel manifolds]{Vector fields on projective Stiefel manifolds\\ and the Browder-Dupont invariant}
\author{Yanghyun Byun}
\address{Department of Mathematics, Hanyang University, Sungdong-gu, Seoul 04763, Korea,}
\email{yhbyun@hanyang.ac.kr}

\author{J\'ulius Korba\v s}\thanks{This work was supported, in case of the second author,
by the Slovak Research and Development Agency under the contract No. APVV-16-0053 and by the grant agency VEGA. 
He was also partially affiliated with the Mathematical Institute, Slovak Academy of Sciences, Bratislava.}
\address{Department of Algebra and Geometry,
Faculty of Mathematics, Physics, and Informatics,
Comenius University, Mlynsk\'a Dolina,
SK-842 48 Bratislava,
Slovakia,}
\email{korbas@fmph.uniba.sk}

\author{Peter Zvengrowski}
\address{Department of Mathematics and Statistics,
The University of Calgary,
Calgary Alberta T2N 1N4,
Canada,}
\email{zvengrow@gmail.com}

\date{}

\begin{abstract}
We develop strong lower bounds for the span of the projective Stiefel manifolds
$X_{n,r}=O(n)/(O(n-r)\times \mathbb Z/2)$, which enable very accurate (in many cases exact) estimates of the span.
The technique, for the most part, involves elementary stability properties of vector bundles.
However, the case $X_{n,2}$ with $n$ odd presents extra difficulties, which are partially resolved using the Browder-Dupont
invariant. In the process, we observe that the symmetric lift due to Sutherland does not necessarily exist for all odd dimensional closed manifolds, and therefore the Browder-Dupont invariant, as he formulated it,
is not defined in general. We will characterize those $n$'s for which
the Browder-Dupont invariant is well-defined on $X_{n,2}$.
Then the invariant will be used in this case
to obtain the lower bounds for the span as a corollary of a stronger result.
\end{abstract}

\maketitle

\section{Introduction, preliminaries, and the main results}
The span of a finite rank real vector bundle $\alpha$, denoted $\text{span}(\alpha)$, is $k$ if
$\alpha$ admits $k$, but no more than $k$, everywhere linearly independent cross sections.
If $\alpha$ is equipped with a Euclidean metric, then $\text{span}(\alpha)\geq k$ means
that $\alpha\approx k\varepsilon\oplus \eta$ for some vector bundle $\eta$;
here and in the sequel $\varepsilon$ is the trivial line bundle and
$k\varepsilon$ denotes the $k$-fold Whitney sum of $\varepsilon$ with itself.

For a $q$-dimensional smooth connected manifold $M^q$, one defines its span to be
$\text{span}(M):=\text{span}(\tau_M)$, where $\tau_M$ is the tangent bundle of $M$.
Sections of $\tau_M$ are usually called vector fields on $M$, and the problem
of determining the number $\text{span}(M)$ is referred to as the vector field problem
on $M$ (\cite{thomas}, \cite{julius-peter}).

Besides the span of a manifold one can consider its stable span (\cite{koschorke},\cite{julius-peter}):
\bdm \text{span}^0(M):=\text{span}(\tau_M\oplus t\varepsilon)-t, \edm
where one may take any $t\geq 1$. Note that $\text{span}^0(M) \geq \text{span}(M)$
always holds.

As we have noticed in \cite{julius-peter}, a possible strategy in solving the vector field problem
on a given manifold is to look for and find (an estimate of) its stable span,
and then to try to show that the same result holds true also for the span.
In the present paper we apply this strategy to obtain strong lower bounds
for the manifolds appearing in the title. In addition to these considerations,
the work will involve computation of the mod 2 Kervaire semi-characteristic
of the projective Stiefel manifolds (see (2.7) below),
as well as the use of the James-Thomas numbers (see Lemma 2.2 and the paragraph above it)
and the Browder-Dupont invariant (see Theorem B, Proposition 2.4 and (3.4) below). Example 1.2 below 
will give the reader a feeling
for the strength of our lower bounds.

The projective Stiefel manifold $X_{n,r}\ (r<n)$ is a closed connected smooth manifold,
obtained from the ordinary Stiefel manifold $V_{n,r}$ of orthogonal $r$-frames in $\mathbb R^n$
by identifying $(v_1,\cdots,v_r)\in V_{n,r}$ with $(-v_1,\cdots,-v_r)$;
we have $\dim(X_{n,r})=\dim(V_{n,r})=nr-\binom{r+1}{2}$. The homeomorphism $X_{n,r}\cong O(n)/(O(n-r)\times \mathbb Z/ 2)$
expresses $X_{n,r}$ as a homogeneous space. In particular, $X_{n,1}$ is the $(n-1)$-dimensional
real projective space $P^{n-1}$, and $X_{n,2}$ can obviously
be interpreted as the total space of the tangent sphere bundle of $P^{n-1}$.

Various topological properties of the projective Stiefel manifolds have been studied 
(cf.~\cite{julius-peter} and
the references there, also e.g.~\cite{barufatti1},\cite{barufatti2}).  
It is well known that  $\text{span}(P^{n-1})=\rho(n)-1$,
where $\rho(n)=2^c+8d$ for $n$ expressed as $(2a+1)2^{c+4d}$, $a, d\geq 0$, $0\leq c \leq 3$ (\cite{adams2}),
hence we confine ourselves to $X_{n,r}$ with $r>1$. The parallelizability question for these
manifolds was almost completely solved in \cite{antoniano2}. Just the parallelizability question for $X_{12,8}$ was left open,
and this remains so up to now. As far as $\text{span}(X_{n,r})$ is concerned,
fairly accurate bounds and even exact values are known in many cases (cf. Example 1.2)
apart from the troublesome case when $n$ is odd and $r=2$, and addressing this case is one of the objectives of this note.

Let $\xi_{n,r}$ be the Hopf line bundle associated to the obvious double covering $V_{n,r}\rightarrow X_{n,r}$.
If no confusion is possible, we write $\xi_{n-1}$ instead of $\xi_{n,1}$. There is also a canonical $(n-r)$-bundle
$\beta_{n,r}$ over $X_{n,r}$, characterized by
\bdm r\xi_{n,r}\oplus \beta_{n,r}\approx n\varepsilon . \tag{1.1}\edm

Writing $\tau_{n,r}$ for the tangent bundle of $X_{n,r}$, we have (\cite{lam}, \cite{peter})
\bdm \tau_{n,r}\approx r\xi_{n,r}\otimes \beta_{n,r}\oplus \binom{r}{2}\varepsilon , \tag{1.2}\edm
and stably
\bdm \tau_{n,r}\oplus \binom{r+1}{2}\varepsilon \approx nr\xi_{n,r}. \tag{1.3}\edm

The the formula (1.2) immediately implies that one always has
\bdm \text{span}(X_{n,r}) \geq \binom{r}{2}. \edm
However, this lower bound usually turns out to be quite weak.

On the other hand, it is quite easy to give a stronger lower bound for stable
$\text{span}(X_{n,r})$.
Indeed, if
\bdm p:X_{n,r}\rightarrow P^{n-1} \edm
is the obvious smooth fibration, then the pullback $p^*(\xi_{n-1})$ is $\xi_{n,r}$, and from (1.3)
we obtain
\bdm \text{span}^0(X_{n,r}) \geq \text{span}(nr\xi_{n-1})-\binom{r+1}{2}:=k_{n,r} .  \tag{1.4}\edm
The span of multiples of the Hopf line bundle $\xi_{n-1}$ over $P^{n-1}$ is well studied,
under the name of generalized vector field problem, and our lower bound $k_{n,r}$ is therefore computable
(see e.g. \cite{lam-1}).

\vskip 8mm
\noindent{\it Example 1.1.} Using \cite[Theorem 1.1]{lam-1}, one finds that $k_{2^m-3, 2}=2^m-3$ for all
$m\geq 3$. On the other hand, we observe that the identity
 $\binom{2^{n+1}-r}{2^n-r}=1 \mod 2$ holds,
when $n\geq 1$ and $0<r\leq 2^n$, which can be proved using Lucas' theorem.
Then (1.3) implies that the Stiefel-Whitney class $w_{2^m-6}(X_{2^m-3,2})$
equals $w_1^{2^m-6}(\xi_{2^m-3,2})$. The latter class does not vanish, because by \cite[Theorem 1.6]{Git-Hand},
the lowest $N$ such that $w_1^N(\xi_{2^m-3,2})=0$ is $N=2^m-4$.
As a consequence, we have
$\text{span}^0(X_{2^m-3,2})=2^m-3$.
\vskip 8mm

Following the above mentioned strategy of relating results for the stable span and results for the span,
one can now naturally ask whether or not $k_{n,r}$ is also a lower bound for $\text{span}(X_{n,r})$.

We let $t$ denote the generator of the ring
$H^*(P^{n-1};\mathbb Z/2)$.

\begin{defn} The mod 2 integer $\chi(n)$ is
such that the $(n-1)$-th Wu class of the stable inverse of $2n\xi_{n-1}$ is $\chi(n)t^{n-1}$.
\end{defn}

The vanishing of $\chi(n)$ for even $n$ is proved in Lemma 5.5, and the table (5.7) below
for $\chi(n)$ with small odd $n$
is also provided. Then one of our main results is as follows.
\vskip 8pt
\noindent{\bf Theorem A.} {\it We have that $\text{span}(X_{n,r})\geq k_{n,r}$ if any of the following conditions hold:
\begin{itemize}[leftmargin=1.5cm]
\item[{\rm (a)}] $r>2$,
\item[{\rm (b)}] $r=2$ and $n$ is even,
\item[{\rm (c)}] $r=2$, $n=3$ or $n$ is odd with $\chi(n)=0$.
\end{itemize} }
\vskip 8pt
This is actually a join of Theorem 2.1, Corollary 2.3, Lemma 2.5 and Theorem B.
We remark that Theorem 2.1 has already been proved in \cite[p.~100]{julius-peter3}, but, since that source
is not generally accessible, we reproduce a short proof in \S 2.
It would also be helpful for the reader to know that $\chi(3)=1$
to see that (c) above is a combination of Corollary 2.3 and Theorem B (see also the table (5.7) below).

The proof of Theorem B is much more complicated.
It is based on applying the Browder-Dupont invariant
in the form defined by W.~Sutherland in \cite{suth}.
We will review Sutherland's construction more closely in \S 3 below
and observe that the condition $\chi(n)=0$ in (c) above
assures that the Browder-Dupont invariant
on $X_{n,2}$ is well-defined.

Here we describe the $(2n-3)$-plane bundle $\sigma_{n,2}$
stably isomorphic to the tangent bundle $\tau_{n,2}$
and such that $\text{span}(\sigma_{n,2})\geq k_{n,2}$. The bundle $2n\xi_{n-1}$ over $P^{n-1}$
can be written $2n\xi_{n-1}\approx \alpha\oplus n\varepsilon$ for a unique (up to isomorphism)
bundle $\alpha$ over $P^{n-1}$, by stability. Let $\sigma_0 = \alpha\oplus (n-3)\varepsilon$
(assuming $n\geq 3$), and let $\sigma_{n,2}=p^*\sigma_0$, a rank $2n-3$ vector bundle over $X_{n,2}$. Note that \ $\sigma_0 \oplus 3\varepsilon
\approx \alpha \oplus n\varepsilon \approx 2n\xi_{n-1}$.
Now $\tau_{n, 2}$ also has rank $2n-3$ and is stably equivalent to $\sigma_{n,2}$ since
\bdm
\sigma_{n,2}\oplus 3\varepsilon \approx p^*(\sigma_0\oplus 3\varepsilon) \approx
p^*(\alpha\oplus n\varepsilon)\approx p^*(2n\xi_{n-1})\approx 2n\xi_{n,2}
\approx \tau_{n,2}\oplus 3\varepsilon. \tag{1.5}
\edm
Since $\text{span}\ \sigma_{n,2} \geq \text{span}(\sigma_0)=-3+\text{span}(2n\xi_{n-1})=k_{n,2}$, it will
suffice to show that $\sigma_{n,2}\approx \tau_{n,2}$ for all odd $n$ with $\chi(n)=0$
to complete the proof of part (c), Theorem A. This is done as follows.

\vskip 8pt
\noindent{\bf Theorem B.} {\it
Assume $n\geq 3$ and $\chi(n)=0$.
Then $b_B(\sigma_{n,2},\rho')=0$. In particular, it follows that $\sigma_{n,2}\approx \tau_{n,2}$.}
\vskip 8pt

\noindent Here $b_B(\sigma_{n,2},\rho')$ denotes the Browder-Dupont invariant (see (3.4) below and Theorem 3.1)
with a suitable trivialization
\bdm \rho':k\varepsilon\stackrel{\approx}{\longrightarrow} \varepsilon\oplus \sigma_{n,2}\oplus \nu_{n,2} \edm
(see (4.2), (5.4), (5.6) below and Remark 4.2)
where $\nu_{n,2}$ is the normal bundle of $X_{n,2}$.
We will also observe that  $b_B(\tau_{n,2}, \kappa)=0$ under the condition $\chi(n)=0$, where $\kappa: k\varepsilon\stackrel{\approx}{\longrightarrow} \tau_{n,2} \oplus \nu_{n,2}$
is the canonical trivialization (see Proposition 2.4). The vanishing of these two invariants proves
that $\sigma_{n,2}\approx \tau_{n,2}$ (see Theorem 3.1 and `Proof of Theorem B' below Lemma 5.4).
Our argument seems to be the first in the literature, in which the Browder-Dupont invariants
are actually calculated
to determine whether or not a pair of two given bundles, each stably isomorphic to the tangent  bundle, are isomorphic to each other.
The calculation occupies about the half of the paper: It begins with the evaluation of $b_B(\tau_{n,2}, \kappa)$ in Proposition 2.4
and occupies the whole of \S 3-5. During the argument, we observe that Sutherland's formulation
of the invariant is in general not well defined (see Condition S in \S 3), which was a surprise
even to the authors themselves.

On the other hand, the following example illustrates the power of the lower bound $k_{n,r}$ in Theorem A.
\vskip 8pt
\noindent {\it Example 1.2.}
\begin{itemize}[leftmargin=1.5cm]
\item[(a)] Taking into account Example 1.1, Theorem A implies that
       $$\text{span}(X_{2^m-3, 2})=k_{2^m-3,2}=2^m-3=\text{span}^0(X_{2^m-3,2})$$
       for all $m\geq 3$ such that $\chi(2^m-3)=0$.
\item[(b)] $\text{span}(X_{14,4})=38$,
\item[(c)] $\text{span}(X_{16,5})=58$,
\item[(d)] $1618\leq \text{span}(X_{58,51})\leq 1625$.
\end{itemize}
\vskip 8pt
We have that $\text{span}(X_{5,2})=5$ and $\text{span}(X_{13,2})=13$ by the example (a) above
and by the table (5.7) below.
In examples (b), (c), (d), the lower bound is furnished in each case by Theorem A.
The upper bounds in (b), (d) are found using Stiefel-Whitney classes, while in (c) the upper bound is
found using the ring structure in {\it K}-theory. A treatment of the upper bounds can be found 
in \cite{julius-peter2} and \cite{sankaran-peter},
which include many infinite families, where the span is exactly found.

We remark that there is an open question, whether or not
$$\text{span}(X_{n,r})=\text{span}^0(X_{n,r})$$
for all projective Stiefel manifolds. The probability that this question will be settled in the affirmative
(and consequently an even stronger result than Theorem A holds) is not low.
In fact, in \cite[Theorem 3.2.11]{julius-peter} a large list is given of those $X_{n,r}$ for which the above equality
can indeed be quite easily shown, using the results of \cite{koschorke}. But approximately one third of the $X_{n,r}$
(in an asymptotic sense) remain undecided.
Therefore Theorem A can be very useful, since it applies to a majority of $X_{n,r}$, including
many of $X_{n,r}$ previously undecided and in particular the $X_{n,2}$ with $n$ odd.

One should not expect that $\text{span}(X_{n,r})=k_{n,r}$ is always true. For instance, by \cite{antoniano2},
the manifold $X_{9,8}$ is parallelizable, so has span $36$, while $k_{9,8}=28$.

The rest of the paper is divided into four sections. In \S2, after having proved Theorem A for the projective Stiefel
manifolds $X_{n,r}$ with $r\geq 3$ as well as for $X_{n,2}$ with $n$ even, we show that from
among the remaining manifolds only $X_{3,2}$ and $X_{5,2}$ have the James-Thomas number 1. The other $X_{n,2}$ have 
James-Thomas number $2$, except possibly for when $n=2^t+1$.
We also prove that the Browder-Dupont invariant
$b_B(\tau_{n,2}, \kappa)$ vanishes under the assumption $\chi(n)=0$.
In \S 3, we briefly review the Browder-Dupont invariant $b_B(\sigma,\rho)$ as formulated by
Sutherland (\cite{suth}).
We will explain why the construction is not applicable
to a general odd-dimensional closed manifold.
In \S 4, we establish an S-duality which will be needed in
the next section.
In \S 5, we prove Theorem B. In particular, we observe that $\chi(n)=0$ when $n$ is even.

We are indebted to Prof.~Parameswaran Sankaran and Prof.~Wilson Sutherland for helpful discussions.
Also many thanks to the two referees for their comments.

%%Section 2%%
\section{A destabilization trick and preparations for applying the Browder-Dupont invariant}
In this section, we prove a major part of Theorem A and make some preparations for
proving its remaining part later in \S 5.

Let $\beta'_{n,r}:=\beta_{n,r}\otimes \xi_{n,r}$. Then the formula (1.1)
and the fact that $\xi\otimes\xi$ is trivial imply that
\bdm \beta'_{n,r}\oplus r\varepsilon\approx n\xi_{n,r}. \tag{2.1}\edm
Further, if $q:X_{n,r}\rightarrow X_{n,r-1}$ is the standard fibration, then obviously
\bdm q^*(\beta_{n,r-1})\approx \beta_{n,r}\oplus\xi_{n,r}, \edm
and therefore
\bdm q^*(\beta'_{n,r-1})\approx \beta'_{n,r}\oplus \varepsilon.\edm
Now suppose $r\geq 3$. From the latter isomorphism we obtain
\bdm q^*(r\beta'_{n,r-1})\approx rq^*(\beta'_{n,r-1})\approx r\beta'_{n,r}\oplus r\varepsilon.\edm
Using the formula (1.2), we then have
\bdm \tau_{n,r}\approx q^*(r \beta'_{n,r-1})\oplus \frac{r(r-3)}{2}\varepsilon. \tag{2.2} \edm
Of course, for some bundle $\alpha$ over $P^{n-1}$ one has by the definition (1.4)
\bdm nr\xi_{n-1} =(\binom{r+1}{2}+k_{n,r})\varepsilon \oplus \alpha. \edm
Multiplying the formula (2.1) by $r-1$, we obtain
\bdm
(r-1)\beta'_{n,r-1}\oplus r(r-1)\varepsilon \approx n(r-1)\xi_{n,r-1}\approx p^*(n(r-1)\xi_{n-1}),
\edm
and hence
\bdm n\beta'_{n,r-1}\oplus r(r-1)\varepsilon \approx
(\binom{r+1}{2}+k_{n,r})\varepsilon\oplus p^*(\alpha). \edm

The latter formula is in dimension $nr$ and can be destabilized to $1+\dim(X_{n,r-1})$.
Now
\bdm nr-\binom{r+1}{2}=\dim (X_{n,r}) \geq 1+\dim (X_{n,r-1}),\edm
hence we can cancel $\binom{r+1}{2}\varepsilon$ and obtain
\bdm r\beta'_{n,r-1}\oplus \frac{r(r-3)}{2}\varepsilon \approx k_{n,r}\varepsilon \oplus p^*(\alpha).\tag{2.3}\edm
Taking $q^*$ of the formula (2.3) and using the formula (2.2), we have
\bdm \tau_{n,r}\approx k_{n,r}\varepsilon\oplus q^*p^*(\alpha).\edm
Hence $\text{span}(X_{n,r})\geq k_{n,r}$ if $r\geq 3$. In addition to this, by
\cite[Theorem 3.2.11]{julius-peter} one has for $n$ even
\bdm \text{span}^0(X_{n,2})=\text{span}(X_{n,2}). \edm
We have thus proved the following part of Theorem A.

%Thm 2.1
\begin{theorem} Let $r>1$. We have $\text{\rm span}(X_{n,r})\geq k_{n,r}$,
except possibly for $X_{n,2}$ with $n$ odd.\end{theorem}

We remark that the case when $n$ is even and $r=2$ in this theorem also admits a proof by the same
destabilization technique as that used for the case $r>2$, taking advantage of the extra fact
that the tangent bundle of $X_{n,1}=P^{n-1}$
admits a nowhere zero section for $n$ even. But the case when $n$ is odd and $r=2$
cannot be handled this way and will be partially
settled in \S 5 by the Browder-Dupont invariant.
This invariant, in addition to the above mentioned methods,
also resolves the case when $n$ is even and $r=2$,
as can be seen by Corollary 5.6.

For an odd dimensional manifold $M^q$ the James-Thomas number is the number of isomorphism classes
of rank $q$
vector bundles over $M$ that are stably isomorphic to the tangent bundle $\tau_M$.
It is shown in \cite{james-thomas} that this number is either $1$ or $2$.
The following lemma determines the James-Thomas number of $X_{n,2}$
for almost all $n$. This will also put our use
of the Browder-Dupont invariant more into context
(see, for instance, the last sentence of \S 5).

%Lemma 2.2
\begin{lemma} The James-Thomas number is $1$ for $X_{3,2}$ and $X_{5,2}$, and it is $2$
for the remaining manifolds $X_{n,2}$, except possibly for $n=2^t+1$, $t\geq 3$.
\end{lemma}
\bp
For $X_{3,2}$ and $X_{5,2}$ it is enough to note that they have dimension $3$ and $7$ respectively, hence
\cite[Theorem 1.7]{james-thomas} gives the result. In the rest of the proof we shall suppose $n\neq 2^t+1$.

Let $BO$ be the classifying space of the stable orthogonal group $O$, and let
$$\sigma: H^{i+1}(BO; \mathbb Z/2)\rightarrow H^i(\Omega BO;\mathbb Z/2)$$
be the suspension homomorphism. In applying \cite[Theorem 1.6]{james-thomas},
we shall replace the loop space $\Omega BO$ by $O$ (see e.g.~\cite[Pretheorem 2.3.1 (iv)]{adams3}).
Then instead of $\sigma(w_{i+1})$, where $w_j$ is the $j$-th universal Stiefel-Whitney class, we shall write
$v_i\in H^i(O;\mathbb Z/2)$. Now it suffices to show that for any map
$\beta:X_{n,2}\rightarrow O$ one has
$$\Delta(\beta):=\beta^*(v_{2n-3})+ \sum_{i=2}^{2n-3}\beta^*(v_{i-1})w_{2n-2-i}(X_{n,2})=0$$
in $\mathbb Z/2$-cohomology.

Any commutative, associative algebra with unit over $\mathbb Z/2$,
generated by elements $x_1, \cdots, x_k$ such that the monomials $x_1^{\varepsilon_1}\cdots x_k^{\varepsilon_k}$,
where $\varepsilon_j=0,1$, form an additive basis, will be written $V(x_1, \cdots, x_k)$.
First observe that by \cite[Theorem 1.6]{Git-Hand}, we have
\bdm H^*(X_{n,2};\mathbb Z/2)=(\mathbb Z/2)[x]/(x^N)\otimes V(x_q), \tag{2.4} \edm
as an algebra, where $x=w_1(\xi_{n,2})$, $N=n-1, n$ according as $n$ is respectively odd, even,
and $q=n-1, n-2$ according as $n$ is respectively odd, even (note $q$ is thus always even).
Since the formula (1.3) implies
\bdm w_{2n-2-i}(X_{n,2})=\binom{2n}{2n-2-i}x^{2n-2-i}, \edm
we have
\bdm \Delta(\beta)=\beta^*(v_{2n-3})+\sum_{i=2}^{2n-3}\beta^*(v_{i-1})\binom{2n}{2+i}x^{2n-2-i}. \edm
Now we recall that by \cite[Theorem 8.7)]{borel} one has for the Steenrod squares
\bdm Sq^i(v_j)=\binom{j}{i}v_{i+j} \edm
for $i\leq j$. Hence it is sufficient to show that
\bdm \beta^*(v_{2^k -1})\in (\mathbb Z/2)[x]/(x^N).\edm
This task is trivial if $k$ is small, and one only has to consider the range for $k$ where $2^k-1>q$.
Then one has
\bdm \beta^*(v_{2^k-1})=\lambda x_q x^{2^k-1-q} \tag{2.5}\edm
for some $\lambda \in \mathbb Z/2$.
By \cite[Theorem 2.1]{antoniano}, one readily checks that $Sq^1(x_q)=0$, and applying $Sq^1$ to the formula (2.5), we obtain
\bdm \lambda x_qx^{2^k-q}= \beta^*(v_{2^k}). \tag{2.6} \edm

Observe that the top class in $H^*(X_{n,2}; \mathbb Z/2)$ is $x_qx^{n-2}$ if $n$ is odd or
$x_qx^{n-1}$ if $n$ is even. Hence in all the cases which we need to consider we have $2^k-q<N$;
therefore $x_q x^{2^k-q} \neq 0$. On the other hand,
\bdm \beta^*(v_{2^k})\in (\mathbb Z/2)[x]/(x^N), \edm
because
\bdm v_{2^k}= Sq^{2^{k-1}}\cdots Sq^2Sq^1(v_1). \edm
Finally, since $2^k\geq q+2 \geq N$, we have $\beta^*(v_{2^k})=0$ and the formula (2.6) implies that $\lambda=0$.
\ep

One can readily check that when an odd dimensional manifold has
the James-Thomas number 1
its stable span and span coincide.
Hence the formula (1.4) and Lemma 2.2 give the following.

%%%%cor2.3
\begin{corollary}
For $X_{3,2}$ and $X_{5,2}$ the span and stable span coincide. In particular,
for $n=3,5$ we have $\text{span}(X_{n,2})\geq k_{n,2}$.
\end{corollary}

On the other hand,  we have $\chi(3)=1$, $\chi(5)=0$ (see the table (5.7))
and, therefore, the case $n=5$ can be proved by Theorem B, \S 1 as well.

We will employ the Browder-Dupont invariant as formulated by Sutherland(\cite{suth}) in order to prove Theorem B. However we will observe that the invariant $b_B(\tau_{n,2}, \kappa)$
is defined only under the condition
$\chi(n)=0$ (see Theorem 3.1 and Lemma 3.2).
Here $\kappa:k\varepsilon\stackrel{\approx}{\longrightarrow} \tau_{n,2}\oplus \nu_{n,2}$ denotes the canonical trivialization, where
$\nu_{n,2}$ is the normal bundle of an embedding $X_{n,2}\rightarrow S^k$ for a $k$ such that
$k\geq 4n-4$.

%%% prop2.4
\begin{proposition}
We have that  $b_B(\tau_{n,2}, \kappa)=0$ under the condition $\chi(n)=0$.
\end{proposition}
\bp
 By  Theorem 3.1 and Lemma 3.2 below (in fact, by \cite[Theorem 2.7]{suth} with a slight modification)
 we know that, under the condition
  $\chi(n)=0$, $b_B(\tau_{n,2}, \kappa)$ is well-defined and is
the Kervaire mod 2 semi-characteristic $\chi_2(X_{n,2})$. We recall that
\bdm \chi_2(X_{n,2})=\sum_{i=0}^{n-2}\dim(H^i(X_{n,2};\mathbb Z/2))\ \ (\text{mod}\ 2).
\tag{2.7}\edm
Hence the proof will be completed by verifying the following result (which is more
general than really needed).
\ep

%%% lemma2.5
\begin{lemma}
Let $r>1$ and $n$ be such that $\dim(X_{n,r})$ is odd. Then one has $\chi_2(X_{n,r})=0$.
\end{lemma}
\bp
We shall distinguish three cases.

(a) Let $n=2k$ and $r=2$. By the formula (2.4) we have now
\bdm H^*(X_{2k,2}; \mathbb Z/2)= (\mathbb Z/2)[x]/(x^{2k})\otimes V(x_{2k-2}).  \edm
hence
\bdm \dim(H^i(X_{2k,2}; \mathbb Z/2))=1 \edm
for $i=0,1,\cdots, 2k-3,$ and
\bdm \dim(H^{2k-2}(X_{2k,2};\mathbb Z/2))=2. \edm
Then the formula (2.7) gives $\chi_2(X_{2k,2})=0$.

(b) Let $n=2k+1$, $r=2$. By the formula (2.4),
\bdm H^*(X_{2k+1, 2}; \mathbb Z/2)=(\mathbb Z/2)[x]/(x^{2k})\otimes V(x_{2k}). \edm
Hence
\bdm \dim(H^i(X_{2k+1,2};\mathbb Z/2))=1 \edm
for $i=0,1,\cdots, 2k-1$, and the formula (2.7) gives the result.

(c) Let $r>2$. We define the following smooth involutions $S$ and $T$ on $X_{n,r}$. If
\bdm w:=\{ (v_1, \cdots. v_r),(-v_1,\cdots, -v_r)\}\in X_{n,r},\edm
then we put
\bdm S(w)=\{(-v_1, v_2,\cdots, v_r), (v_1,-v_2, \cdots, -v_r)\} \edm
and
\bdm T(w)=\{(v_1,-v_2, v_3, \cdots,v_r), (-v_1,v_2,-v_3, \cdots, -v_r)\}. \edm
Then $ST(=TS)$, $S$ and $T$ are fixed point free involutions on $X_{n,r}$, and
together with the identity map they give a free action of the group
$\mathbb Z/2 \times \mathbb Z/2$ on $X_{n,r}$.
As a consequence of R.~Stong's lemma \cite[Lemma 3.4]{stong},
we have then $\chi_2(X_{n,r})=0$.
\ep

%%%3
\section{The Browder-Dupont invariant and a review of Sutherland's symmetric lift}
We recall the Browder-Dupont invariant
as formulated by W.~Sutherland (\cite{suth}).
In general, let $M$ be a connected closed smooth manifold with an odd dimension $q$.
We may assume $M\subset S^k$ with $k\geq 2q+2$.
Let $\tau$ and $\nu$ denote respectively the tangent and normal bundle of this embedding.
We consider the universal real vector bundle $\bar\gamma$ of rank $2k-2q$
whose $(q+1)$'st Wu class vanishes.
We choose also a fixed S-duality between a space $X$ and (a finite skeleton of) $T(\bar\gamma)$.

There is
the normal invariant $S^k\rightarrow T(\nu)$, which is realized in this paper
as the collapse map $ S^k \rightarrow N/\partial N \cong T(\nu)$,
where $N$ is the tubular neighborhood of $M$ in $S^k$.
Furthermore assume
$\sigma$ and $\eta$ are vector bundles over $M$ such that
$l\varepsilon\approx \sigma\oplus\eta$ for some integer $l>0$.
Then $T(\sigma)$
and $T(\eta\oplus \nu)$ are S-dual to each other where the S-duality is given
by the composite (e.g.~\cite[above Proposition 3.1]{dupont}),
\bdm
\begin{aligned}
S^{l+k}\equiv \Sigma^l S^k & \rightarrow \Sigma^l T(\nu)\equiv T(\varepsilon^l\oplus \nu)\\
&\cong T(\sigma\oplus\eta \oplus \nu)
\rightarrow T(\sigma)\wedge T(\eta\oplus\nu).
\end{aligned}
\tag{3.1}
\edm

On the other hand,  Sutherland, by Proposition 2.1 in \cite{suth}, considers
a {\it symmetric} lift $a:\nu\times \nu\rightarrow \bar\gamma$
of the classifying map $\nu\times\nu\rightarrow \gamma$, where
$\gamma$ is the universal real vector bundle.
If $t$ is the bundle map from
$\nu\times \nu$ to itself transposing the factors,
then $a$ is symmetric in the sense that
$a$ and $a\circ t$ are homotopic to each other through bundle maps.
Consider the natural bundle map $\bar\Delta:2\nu\rightarrow \nu\times\nu$
covering the diagonal map $M\rightarrow M\times M$.
Then he shows that the bundle map $a\circ\bar\Delta: 2\nu \rightarrow \bar\gamma$
does not depend on the choice of the symmetric lift up to homotopy.
In fact most of these statements can be found in a more general form
in the paragraph before
\cite[Proposition 2.1]{suth}.

However, it does not appear that the symmetric lift exists for all odd-dimensional
connected closed manifolds. After a closer look at \cite{suth} one notices that it exists under the following condition:
\begin{quote}
{\bf Condition S}: The $i$-th Wu class $v_i(\nu)$ vanishes
for each $i$ such that $\frac{q+1}{2} \leq i \leq q$.
\end{quote}
In fact, the above condition can be seen by combining the first sentence
of the last paragraph of \cite[\S 4]{suth} with the first two lines of the last paragraph of the previous page.
In spite of what is said in \cite{suth}, Condition S above is an extra condition imposed on $M$, as shown
by the following example.
If $\tau$ and $\nu$ respectively denote
the tangent bundle and the normal bundle of
the real projective space $P^5$, regarding the total Wu classes, $v(\tau)$ and $v(\nu)$,
we have that respectively
$v(\tau)=1+t^2$ and $v(\nu)=1+t^2 +t^4$. Here $t\in H^1(P^5; \mathbb Z/2)$
denotes the generator of $H^*(P^5;\mathbb Z/2)$. Therefore Condition S
is not satisfied by $P^5$. The origin of the assertion that Condition S holds for general $M$ might be
the fact that $v_i(\tau)=0$ for the same range of $i$'s.
The same mistake is present in Dupont's work (\cite{dupont}, in particular, 
the 9th line from the bottom on p.~219), which might be the true
origin of the mistake of \cite{suth}. So far, in spite of our efforts,
we have not been able to recover the general case. For instance,
we considered employing the stable inverse of the universal bundle,
but to no avail. 

Let $\sigma$ be a rank $q$ vector bundle stably isomorphic to $\tau$. Then there
is a trivialization
\bdm \rho: k\varepsilon \rightarrow \sigma\oplus\nu. \tag{3.2}\edm
Consider two pairs, $(\sigma_1,\rho_1)$ and $(\sigma_2,\rho_2)$,
 where $\sigma_i$, $i=1,2$, are rank $q$ bundles stably isomorphic to $\tau$ and
$\rho_i:k\varepsilon \rightarrow \sigma_i\oplus\nu$, $i=1,2$,
are trivializations.
Then one says that the two pairs are {\it equivalent} if there is an isomorphism 
$\varphi:\sigma_1\rightarrow\sigma_2$ such that
$(\varphi\oplus 1)\rho_1$ is homotopic to $\rho_2$ through bundle isomorphisms.

In fact, Sutherland uses a convention, which specifies
a bundle equivalence
$$\alpha:\sigma \oplus 2\varepsilon \rightarrow \tau \oplus 2\varepsilon$$
and applies it to the canonical trivialization of $\tau\oplus\nu$
to specify a trivialization of $\sigma\oplus\nu$
(\cite[\S 2]{suth}). Up to homotopy
through isomorphisms
any trivialization of
$\sigma \oplus \nu$ can be obtained in this way (see \cite[Lemma 2.8]{dupont})
since the rank of $\nu$ is larger than $q+1$ by our assumption.

Assume Condition S.
Then there is a symmetric lift $a:\nu\times\nu \rightarrow
\bar\gamma$.
Note that the normal invariant, $S^k\rightarrow T(\nu)$,
together with a trivialization $\rho:k\varepsilon\rightarrow \sigma\oplus\nu$, determines an S-duality
$S^{2k}\rightarrow T(\sigma)\wedge T(2\nu)$ by (3.1) above.
This duality determines the dual $g:X\rightarrow \Sigma^mT(\sigma)$ of $T(a\circ \bar\Delta):T(2\nu)\rightarrow T(\bar\gamma)$, for
an appropriate integer $m$.
Write $f$ for the composite,
\bdm
\begin{CD}
X @>g>>\Sigma^mT(\sigma)@>{\Sigma^m U_{\sigma}}>> \Sigma^m K_q,\end{CD} \tag{3.3}
\edm
where $K_q$ is the Eilenberg-MacLane space of type $(\mathbb Z/2,q)$ and $U_{\sigma}$ is the Thom class.
It should be understood here and in the sequel
that all cohomology is taken with
$\mathbb Z/2$-coefficients.
Then one defines the Browder-Dupont invariant by
 \bdm b_B(\sigma,\rho)=Sq^{q+1}_f(\Sigma^m\iota)\in H^{m+2q}(X)\cong\mathbb Z/2. \tag{3.4} \edm
In the above, $Sq^{q+1}_f(\Sigma^m\iota)$ is the functional Steenrod  square,
 defined by W.~Browder \cite{browder}, based on the work by N.~Steenrod \cite{steenrod},
and $\iota\in H^q(K_q)$ denotes the characteristic element.
Then Sutherland has proved the following
(see \cite[Theorems 2.3, 2.5 and 2.7]{suth}).

%Theorem 3.1
\begin{theorem} {\rm (Sutherland)}
Assume that a symmetric lift $\nu\times\nu\rightarrow \bar\gamma$ exists. Then
 $b_B(\sigma,\rho)$ is well defined for any pair $(\sigma, \rho)$ and
 the equality $b_B(\sigma_1,\rho_1)=b_B(\sigma_2,\rho_2)$ holds if and only if
the pairs $(\sigma_1,\rho_1)$, $(\sigma_2,\rho_2)$ are equivalent. In addition,
the equality, $b_B(\tau, \kappa)=\chi_2(M)$, holds  for the tangent bundle $\tau$ of $M$
and the canonical trivialization $\kappa:k\varepsilon\rightarrow \tau\oplus\nu$.
\end{theorem}
In particular, the equality $b_B(\sigma_1,\rho_1)=b_B(\sigma_2,\rho_2)$ implies that
$\sigma_1\approx\sigma_2$. Also note that $\chi_2(M)$ denotes the $\text{mod}\ 2$ semi-characteristic
of $M$ (see the formula (2.7)).

Recall $\chi(n)\in \mathbb Z/2$ given by Definition 1.1.
Let $\nu_{n,2}$ represent the stable normal bundle of $X_{n,2}$.

%Lemma 3.2
\begin{lemma}
If $\chi(n)=0$, then there is a symmetric lift,
$\nu_{n,2}\times\nu_{n,2} \rightarrow \bar\gamma$.
\end{lemma}
\bp
Note that $\tau_{n,2}$ is stably
isomorphic to $p^*(2n\xi_{n-1})$ and thus $\nu_{n,2}$ is stably isomorphic
to the pull-back of a bundle on $P^{n-1}$.
Hence every characteristic class of $\nu_{n,2}$ lives in $p^*H^*(P^{n-1})$.
Therefore, if $\chi(n)=0$, then Condition S above, that is, that $v_i(\nu_{n,2})=0$
for $i\geq \frac{(2n-3)+1}{2}=n-1$ is satisfied and a symmetric lift exists.
\ep

It could be true that the Browder-Dupont invariant as formulated
in \cite{suth}
may have been used incorrectly on some occasions. For instance,
it is used to define the tangent fibration of a Poincar\'e complex of odd dimension in
\cite{byun1} by one of the present authors,
which however now appears to lack full generality.
It has been
also used in \cite{byun2} to prove the existence of a homotopy equivalence
$M_1\rightarrow M_2$ between closed smooth manifolds
which is not covered by any bundle map
between the tangent bundles. In this case the argument remains valid
if one chooses and fixes a lift $2\nu\rightarrow \bar\gamma$.
In any case the invariant, as formulated by Sutherland in \cite{suth},
still appears useful as illustrated by Theorems A and B above.

%%%% section 4 %%%
\section{An S-duality}
In this section we assume $M_1^n \subset M_2^{n+l}$ are smooth, connected and closed manifolds
with $l>0$.
We let $\tau_1$ and $\tau_2$
denote the respective tangent bundles. Furthermore we assume the normal bundle
of $M_1$ in $M_2$ is trivial. Then it follows that
$\tau_2|_{M_1}\approx l\varepsilon\oplus \tau_1$.
We regard $M_2$ as a submanifold of $S^k$ with $k\geq 2n+2l+2$ and
let $\nu_2$ denote the normal bundle. Also we let $\nu_1:=\nu_2|_{M_1}$.
Then $l\varepsilon \oplus\nu_1$ can be identified with the normal bundle of $M_1$ in $S^k$.

Let $\sigma_2$ be a bundle of rank $n+l$, stably isomorphic to $\tau_2$.
Then there is a trivialization,
\bdm \rho: k\varepsilon \rightarrow \sigma_2 \oplus \nu_2 . \tag{4.1}\edm
Also let $\sigma_1$ be a rank $n$ bundle over $M_1$ such that
there is an isomorphism $\sigma_2|_{M_1}\rightarrow l\varepsilon\oplus \sigma_1$,
which implies that $\sigma_1$ is stably isomorphic to $\tau_1$.
It follows that there is a restriction of $\rho$,
\bdm \rho':k\varepsilon\rightarrow l\varepsilon \oplus \sigma_1 \oplus \nu_1 ,\tag{4.2}\edm
which, in fact, depends on the choice of an isomorphism $\sigma_2|_{M_1}\rightarrow l\varepsilon \oplus \sigma_1$.
By applying (3.1) these two trivializations respectively determine two S-dualities.
One of the two is
\bdm S^{2k}\rightarrow T(\sigma_2)\wedge T(2\nu_2).\tag{4.3} \edm
The other is
\bdm S^{2k}\rightarrow (\Sigma^l T(\sigma_1))\wedge T(2\nu_1).\tag{4.4}\edm
For this duality, we let $N\subset M_2$ denote a tubular neighborhood of $M_1$.
Then there is a homeomorphism,
\bdm T(\nu_2|_N)/T(\nu_2|_{\partial N})\cong
T(l\varepsilon\oplus \nu_1)\equiv \Sigma^l T(\nu_1) . \tag{4.5}\edm
The normal invariant $S^k \rightarrow \Sigma^lT(\nu_1)$, which determines
the S-duality (4.4), is the composite of the two collapses,
$$S^k \longrightarrow T(\nu_2) \rightarrow T(\nu_2|_N)/T(\nu_2|_{\partial N}),$$
followed by the homeomorphism (4.5).

In general, we have a homeomorphism,
\bdm T(\zeta|_N)/T(\zeta|_{\partial N})\cong T(l\varepsilon \oplus \zeta|_{M_1})\equiv \Sigma^l T(\zeta|_{M_1})\edm
for any vector bundle $\zeta$ over $M_2$, in the same way as (4.5).
Here we specify this map as follows: First of all, we choose a homeomorphism $N\to D^l\times M_1$.
Then a retraction $r:N\rightarrow M_1$ is determined. Also, by choosing an isomorphism
$\zeta|_N\to r^*(\zeta|_{M_1})$, a bundle map $r_*:\zeta|_N\rightarrow\zeta|_{M_1}$
covering $r$ has been determined. Then a homeomorphism $D(\zeta|N)\to D^l\times D(\zeta|_{M_1})$
is determined by $r_*$ and the map $N\to D^l\times M_1$, where $D(\cdot)$ denotes the disk bundle.
Note that the homeomorphism $D(\zeta|N)\to D^l\times D(\zeta|_{M_1})$ maps $S(\zeta|N)\cup D(\zeta|_{M_1})$
to $S^{l-1}\times D(\zeta|_{M_1})$ $\cup D^l\times S(\zeta|_{M_1})$, where $S(\cdot)$ denotes
the sphere bundle. This gives the homeomorphism above.

Since there is the collapse map $T(\zeta)\to T(\zeta)/T(\zeta|_{M_2-\text{int}N})$ $\equiv T(\zeta|_N)/T(\zeta|_{\partial N})$
we have the collapse map, $T(\zeta)\to \Sigma^lT(\zeta|_{M_1})$.
Let $c:T(\sigma_2)\rightarrow \Sigma^lT(\sigma_1)$
denote the collapse map.
Let $i_*:\nu_1\rightarrow \nu_2$ and $2i_*:2\nu_1\rightarrow 2\nu_2$ denote the inclusions.
Then the main result of the section is as follows.

%lemma 4.1
\begin{lemma}
The collapse map,
$c:T(\sigma_2)\rightarrow \Sigma^lT(\sigma_1)$, is S-dual to the map,
$T(2i_*): T(2\nu_1)\rightarrow T(2\nu_2)$, with respect to the dualities,
(4.3) and (4.4).
\end{lemma}
\bp
It is enough to show that the diagram below commutes up to homotopy:
\bdm
\begin{CD}
S^{2k} @>>>   T(\sigma_2)\wedge T(2\nu_2)\\
         @VVV       @Vc\wedge 1VV  \\
(\Sigma^{l}T(\sigma_1))\wedge T(2\nu_1) @>1\wedge T(2i_*)>>
(\Sigma^{l}T(\sigma_1))\wedge  T(2\nu_2) ,
\end{CD}
\tag{4.6}
\edm
where the top row and
the left column are respectively the duality maps given by (4.3) and
by (4.4).

We begin by observing that the diagram below commutes:
\bdm
\begin{CD}
T(\sigma_2 \oplus 2\nu_2) @>>> T(\sigma_2)\wedge T(2\nu_2)\\
        @V{\bar c}VV                    @Vc'\wedge 1 VV \\
T((\sigma_2 \oplus 2\nu_2)|_N)/T((\sigma_2 \oplus 2\nu_2)|_{\partial N})
                            @>>>  (T(\sigma_2|_N)/T(\sigma_2|_{\partial N}))\wedge T(2\nu_2)\\
        @V=VV          @A{1\wedge T(2\bar i_*)}AA \\
T(\sigma_2\oplus 2\nu_2)|_N)/T((\sigma_2 \oplus 2\nu_2)|_{\partial N})
                           @>>> (T(\sigma_2|_N)/T(\sigma_2|_{\partial N}))\wedge T(2\nu_2|_N).
\end{CD}
\tag{4.7}
\edm
In the above, $\bar c$ and $c'$ are the obvious collapses
and $2\bar i_*:2\nu_2|_N\rightarrow 2\nu_2$ is the inclusion.
The rows represent the maps induced
by diagonal maps (see, for instance, \cite[p.~207]{dupont}).
Here, by a diagonal map, we mean a restriction of the diagonal. For example,
the middle row is induced by the restriction $N\rightarrow N\times M_2$
of the diagonal $M_2\rightarrow M_2\times M_2$.

Now let $\lambda$ denote the homeomorphism, $T(\sigma_2|_N)/T(\sigma_2|_{\partial N})\rightarrow \Sigma^lT(\sigma_1)$.
We consider the diagram:
\bdm
\begin{CD}
T((\sigma_2\oplus 2\nu_2)|_N)/T((\sigma_2\oplus 2\nu_2)|_{\partial N})
                                      @>>> (T(\sigma_2|_N)/T(\sigma_2|_{\partial N}))\wedge T(2\nu_2|_N)\\
    @V\cong VV                       @A\lambda^{-1}\wedge T(i'_*)AA        \\
\Sigma^lT(\sigma_1\oplus 2\nu_1) @>>>(\Sigma^lT(\sigma_1))\wedge T(2\nu_1) .
\end{CD}
\tag{4.8}
\edm
Here $i'_*:2\nu_1 \rightarrow 2\nu_2|_N$ means the inclusion.
This diagram cannot commute since the upper row is induced by
the diagonal $N\rightarrow N\times N$ while $i'_*$ in the right column covers the inclusion
$M_1 \hookrightarrow N$.  However it commutes exactly
if we replace $\lambda^{-1}\wedge T(i'_*)$
with the map,
$$\lambda \wedge T(r_*):(T(\sigma_2|_N)/T(\sigma_2|_{\partial N}))\wedge T(2\nu_2|_N)
\rightarrow (\Sigma^lT(\sigma_1))\wedge T(2\nu_1),$$
where $r_*:(2\nu_2)|_N\rightarrow 2\nu_1$ is the bundle map covering a retraction $N\cong D^l\times M_1\rightarrow M_1$.
In fact we may write $[a, x, v_x, w_x]$ to denote a point of
$T((\sigma_2\oplus 2\nu_2)|_N)/T((\sigma_2\oplus 2\nu_2)|_{\partial N})$, where $(a, x) \in D^l \times M_1 \cong N$ and $(v_x, w_x)$ is an element of the fiber
$D(\sigma_1\oplus 2\nu_1)_x \equiv D(\sigma_1)_x\times D(2\nu_1)_x$ of the disk bundle. This point is mapped to
$[a, x, v_x, w_x]\in (\Sigma^lT(\sigma_1))\wedge T(2\nu_1)$ by both of the composites in the new diagram.
Also note that $T(r_*)$ and $T(i'_*)$ are homotopy inverses to each other. Thus (4.8) commutes
up to homotopy.

The upper row of (4.8) is the bottom row of (4.7).
Therefore we may combine (4.8) with (4.7) to have a diagram which commutes up to homotopy.
Furthermore, we replace  $T(\sigma_2|_N)/T(\sigma_2|_{\partial N})$
with $\Sigma^l T(\sigma_1)$ exploiting the homeomorphism
$\lambda:T(\sigma_2|_N)/T(\sigma_2|_{\partial N})\rightarrow \Sigma^l T(\sigma_1)$.
Also we write $T((\sigma_2 \oplus 2\nu_2)|_N)/T((\sigma_2\oplus 2\nu_2)|_{\partial N})=
\overline{T((\sigma_2 \oplus 2\nu_2)|_N)}$. Then
(4.7) becomes the following diagram, which commutes up to homotopy:
\bdm
\begin{CD}
T(\sigma_2 \oplus 2\nu_2) @>\Delta_2 >> T(\sigma_2)\wedge T(2\nu_2)\\
         @V\bar c VV               @Vc\wedge 1 VV\\
\overline{T((\sigma_2 \oplus 2\nu_2)|_N)}
                            @>\Delta_N>>  (\Sigma^lT(\sigma_1))\wedge T(2\nu_2)\\
          @V\cong VV               @A{1\wedge T(2i_*)}AA\\
\Sigma^l T(\sigma_1\oplus 2\nu_1) @>\Delta_1>> (\Sigma^lT(\sigma_1))\wedge T(2\nu_1),
\end{CD}
\edm
where $\Delta_2$, $\Delta_N$ and $\Delta_1$ originate from the maps induced by diagonal maps.
In particular, we note that $\Delta_N$ cannot be induced by diagonal map
in the same sense as in the second sentence below diagram (4.7).

Let $\hat c:T(\nu_2)\rightarrow T(\nu_2|_N)/T(\nu_2|_{\partial N})$
be the obvious collapse.
Write  $\overline{T((k\varepsilon\oplus \nu_2)|_N)}$
to denote the quotient space $T((k\varepsilon\oplus \nu_2)|_N)/T((k\varepsilon\oplus \nu_2)|_{\partial N})$.
Recall the trivializations $\rho: k\varepsilon \rightarrow \sigma_2 \oplus \nu_2$
and $\rho':k\varepsilon\rightarrow l\varepsilon \oplus \sigma_1 \oplus \nu_1$
respectively from (4.1)
and from (4.2). Then consider the bundle maps, $\rho\oplus 1: k\varepsilon\oplus \nu_2 \rightarrow \sigma_2 \oplus \nu_2\oplus\nu_2$ and
$\rho'\oplus 1:k\varepsilon\oplus\nu_1\rightarrow l\varepsilon \oplus \sigma_1 \oplus \nu_1\oplus \nu_1$.
The map $\overline{T(\rho\oplus 1)}$ is induced by restricting $\rho\oplus 1$ as a bundle map from
$(k\varepsilon\oplus \nu_2)|_N$ to $(\sigma_2\oplus 2\nu_2)|_N$.
Recall
the normal invariants
$S^k\rightarrow \Sigma^lT(\nu_1)$
and $S^k\rightarrow T(\nu_2)$ which induces the S-dualities (4.4)
and (4.3). We denote these normal invariants respectively by $\varphi_1$ and $\varphi_2$.
Then it is straightforward to see that
the following diagram commutes up to homotopy:
\bdm
\begin{CD}
S^{2k}@>\Delta_2\circ T(\rho\oplus 1)\circ\Sigma^k\varphi_2>> T(\sigma_2)\wedge T(2\nu_2)\\
@V=VV                                       @V c\wedge 1 VV\\
S^{2k} @>\Delta_N\circ \overline{T(\rho\oplus 1)}\circ \Sigma^k(\hat c \circ \varphi_2)>> (\Sigma^lT(\sigma_1))\wedge T(2\nu_2)\\
@V=VV                                       @A{1\wedge T(2i_*)}AA\\
 S^{2k} @>\Delta_1\circ T(\rho'\oplus 1)\circ \Sigma^k\varphi_1>>
                              (\Sigma^lT(\sigma_1))\wedge T(2\nu_1).\\
\end{CD}
\edm

It follows that the diagram (4.6) commutes up to homotopy.
\ep

\noindent {\bf Remark 4.2.} The Browder-Dupont invariant depends on
the trivialization $k\epsilon\rightarrow \sigma\oplus \nu$
(see (3.2) and (3.4) above), in particular, when the James-Thomas number of the odd-dimensional manifold in concern is 1
(see {\it Proof of Theorem 2.5} in \cite{suth}). In our case, the trivialization
$\rho': k\epsilon\rightarrow l\epsilon\oplus\sigma_1\oplus \nu_1
\equiv\sigma_2|_{M_1}\oplus \nu_2|_{M_1}$ is chosen as the restriction
of $\rho: k\epsilon \rightarrow \sigma_2\oplus\nu_2$ (see (4,1), (4.2) above and (5.6) below).
In other words, what is special about $\rho'$ in the current manuscript is that it can be extended to a bundle map $k\epsilon \rightarrow \sigma_2\oplus\nu_2$.

%%% section 5
\section{The calculation of the Browder-Dupont invariant}
For each $n\geq 3$, we introduce a manifold,
\bdm Z_n=(S^{n-1}\times S^{n-1})/\{\pm 1\}. \edm
Let $[x,y]$ denote $\{ \pm(x,y) \}\in Z_n$ for any $(x,y)\in S^{n-1}\times S^{n-1}$.
The tangent bundle $\bar\tau_n$ of $Z_n$ may be identified with the set consisting of
$\{\pm(x,y; v,w)\}$ where
$(x,y)\in S^{n-1}\times S^{n-1}$ and $v, w \in \mathbb R^n$ are such that
$x\cdot v =0$ and $y\cdot w=0$.
We also consider the iterated sum
$2n\varepsilon = (n\varepsilon)\oplus (n\varepsilon)$ of the trivial line bundle over $Z_n$.
Then let $\eta_1$ and $\eta_2$ be the line subbundles of $2n\varepsilon$
which respectively have the
fibers $\{[x,y]\}\times \langle x,0 \rangle$
and $\{[x,y]\}\times \langle 0,y\rangle$ at any $[x,y]\in Z_n$. Here $\langle a, b \rangle$
denotes the one dimensional subspace of $\mathbb R^n \oplus \mathbb R^n$ generated by a non-zero $(a,b)\in \mathbb R^n \oplus \mathbb R^n$.
We denote by  $\bar\eta_Z$ the subbundle $\eta_1 +\eta_2$ of $2n\varepsilon$, and by $\beta_Z$, its orthogonal complement.
Let $\bar \xi_n$ denote the Hopf line bundle of $Z_n$.
Then we have, in a way similar to \cite[4.4]{milnor}, that
$${\rm Hom}(\bar\xi_n, \beta_Z)\equiv \bar\tau_n .$$
Also note that $\bar\xi_n$ is canonically isomorphic to $\eta_i$,
for each $i=1, 2$. Therefore, we have that
\bdm
2n\bar\xi_n
\approx {\rm Hom}(\bar\xi_n, 2n\varepsilon)\approx {\rm Hom}(\bar\xi_n, \beta_Z\oplus \bar\eta_Z )\approx \bar\tau_n\oplus 2\varepsilon.
\tag{5.1}
\edm

Let $\bar p$ be the projection $\bar p:Z_n\rightarrow P^{n-1}$ defined by $\bar p[x,y]=[x]$
for any $[x,y]\in Z_n$ and recall the canonical line bundle $\xi_{n-1}$ over $P^{n-1}$.
Then we have that $\bar p^*\xi_{n-1}\approx \bar\xi_n$.
Furthermore consider the inclusion $i: X_{n,2}\hookrightarrow Z_n$,
and recall the projection $p:X_{n,2}\rightarrow P^{n-1}$
(see just above (1.4)).
Then we have that $2n\xi_{n,2}\approx i^*(2n \bar\xi_n)$ since $\bar p i=p$.
Also recall that $2n\xi_{n,2}$ is stably isomorphic to $\tau_{n,2}$ (see (1.3)).
Then considering (5.1),
we have that both $\tau_{n,2}$ and $i^*\bar\tau_n$ are stably isomorphic
to $2n \xi_{n,2}$.

In particular, both $\tau_{n,2}$ and $\bar \tau_n$ are stably pull-backs of
$2n\xi_{n-1}$ and therefore are orientable. Thus the normal bundle
of $X_{n,2}$ in $Z_n$ is orientable and we have:
%5.1
\begin{lemma}
The normal bundle of $X_{n,2}$ in $Z_n$ is a trivial line bundle.
\end{lemma}

We may assume that $Z_n$ is embedded in a sphere $S^k$ where $k \geq 4n-2$
and that $\bar\nu_n$ is the normal bundle.
Let $N\subset Z_n$ denote a tubular neighborhood of $X_{n,2}$.
Then, there is a homeomorphism:
\bdm T(\bar\nu_n|_N)/T(\bar\nu_n|_{\partial N})\cong
T(\varepsilon\oplus i^*\bar\nu_n)\equiv \Sigma T(i^*\bar\nu_n) . \tag{5.2}\edm
Let us denote $i^*\bar\nu_n$ by $\nu_{n,2}$ which represents
the normal bundle of $X_{n,2}$.

We recall from (1.5) the rank $2n-3$ subbundle $\sigma_0$ of $2n\xi_{n-1}$  such that
$2n\xi_{n-1}\approx \sigma_0\oplus 3\varepsilon$.
Then the vector bundle
$\sigma_{n,2}$ over $X_{n,2}$, which is defined as $p^*\sigma_0$,
is
stably isomorphic to $2n\xi_{n,2}$.
Now we introduce a vector bundle
$\bar \sigma_n = \bar p^*\sigma_0$ over $Z_n$.
Also note that $\bar\xi_n\approx \bar p^*\xi_{n-1}$
which means $\bar\sigma_n=\bar p^*\sigma_0$ is stably isomorphic
to $2n \bar p^*\xi_{n-1}\approx  2n\bar\xi_n$.
In addition $\bar\tau_n$ is stably isomorphic to $2n\bar\xi_n$ by (5.1) above.
Thus both $\bar\sigma_n$ and $\bar\tau_n$ are stably isomorphic to $2n\bar\xi_n$.
Therefore $\bar\sigma_n \oplus \bar\nu_n$
is a trivial bundle of rank $k-1$.
Therefore, there is an isomorphism,
\bdm \rho: k\varepsilon\rightarrow \varepsilon\oplus \bar\sigma_n\oplus\bar\nu_n. \tag{5.3}\edm
Since we have $\sigma_{n,2}=i^*\bar\sigma_n$ and $\nu_{n,2}=i^*\bar\nu_n$,
$\rho$ restricts to
\bdm \rho':k\varepsilon\rightarrow \varepsilon\oplus \sigma_{n,2}\oplus \nu_{n,2}.\tag{5.4} \edm

Let $2i_*:2\nu_{n,2}\rightarrow 2\bar\nu_n$ denote the inclusion
and $c:T(\bar\sigma_n)\rightarrow \Sigma T(\sigma_{n,2})$ be the collapse map.
Given the trivializations (5.3) and (5.4), we may consider the S-dualities (4.3) and (4.4) above to apply Lemma 4.1.
Then we have the following.

%lemma 5.2
\begin{lemma} The dual of the map $T(2i_*): T(2\nu_{n,2})\hookrightarrow T(2\bar\nu_n)$
is the collapse $c:T(\bar\sigma_n)\rightarrow \Sigma T(\sigma_{n,2})$.
\end{lemma}

In general, suppose that $\bar\zeta_n$ is a bundle over $Z_n$, of an arbitrary rank $m$, stably isomorphic to $2n\bar\xi_n$.
Write $\zeta_{n,2}=i^*\bar\zeta_n$.
Then there is the collapse map,
$$c:T(\bar\zeta_n)\rightarrow \Sigma T(\zeta_{n,2}).$$
Let $U_{n,2} \in H^m(T(\zeta_{n,2}))$ be the Thom class. Then we observe the following.

%lemma 5.3
\begin{lemma} Assume $n\geq 3$. Then we have that
$$c^*(\Sigma U_{n,2})=0\in H^{m+1}(T(\bar\zeta_n)).$$ \end{lemma}
\bp
Note that
$$c^*(\Sigma U_{n,2})\in H^{m+1}(T(\bar\zeta_n))=\{0, U_Z\cup a\}$$
where $U_Z\in H^m(T(\bar\zeta_n))$ is the Thom class and $a$ is the generator of $H^1(Z_n)$.

Now consider the diagonal map $\Delta: P^{n-1}\rightarrow Z_n$ defined by
$\Delta[x] = [x,x]$. Also recall the projection $\bar p:Z_n\rightarrow P^{n-1}$ given by $\bar p[x,y]=[x]$.
Then $\bar p\circ \Delta$ is the identity on $P^{n-1}$,
which means that $\bar p^*:H^*(P^{n-1})\rightarrow H^*(Z_n)$ is injective.
Note that $a=\bar p^*t$
where $t$ is the generator of $H^*(P^{n-1})$.
Thus we have that
$a^2\neq 0$.

The first Stiefel-Whitney class $w_1$ of $\zeta_{n,2}$
is zero since it is stably isomorphic to $2n\xi_{n,2}$.
Therefore we have that $Sq^1 U_{n,2}= U_{n,2}\cup w_1=0$, and conclude that $Sq^1 c^*(\Sigma U_{n,2})=0$.
On the other hand, we have  $Sq^1 (U_Z\cup a)=U_Z\cup a^2\neq 0$.
Thus $c^*(\Sigma U_{n,2})= U_Z\cup a$ is impossible. Therefore we have that
$c^*(\Sigma U_{n,2})=0$.
\ep

Now let $K_m$ denote the Eilenberg-MacLane space $K(\mathbb Z/2, m)$ for any positive integer $m$.
The lemma above tells us that the composite
\bdm
\CD T(\bar\sigma_n)@>c>>\Sigma T(\sigma_{n,2})@>\Sigma U_{n,2}>> K_{2n-2}
\endCD
\edm
is homotopic to the constant.
In fact we have the following.

%5.4
\begin{lemma}
Assume $n\geq 3$. Then the composite,
\bdm
\CD
T(\bar\sigma_n)@>c>>  \Sigma T(\sigma_{n,2}) @>\Sigma U_{n,2}>> \Sigma K_{2n-3},
\endCD
\edm
is homotopic to the constant.
\end{lemma}
\bp
Recall the notation of (1.5) above.
Then, for a dimensional reason, we have $\sigma_0 \approx \beta\oplus (n-2)\varepsilon$ over $P^{n-1}$,
for some rank $n-1$ vector bundle $\beta$,
which is not necessarily unique.
Consider the vector bundles $\bar\sigma_n'=\bar p^*(\beta\oplus (n-3)\varepsilon)$ and $\sigma_{n,2}'=i^*\bar\sigma_n'$, which are of rank $2n-4$.
Then we have that
\bdm
\varepsilon\oplus \bar\sigma_n' \approx \bar\sigma_n  {\rm\ \ and\ \ }
\varepsilon\oplus \sigma_{n,2}' \approx \sigma_{n,2}, \tag{5.5}
\edm
where we choose the second isomorphism as the restriction of the first.
We have the collapse map,
$$c':T(\bar\sigma_n')\rightarrow \Sigma T(\sigma_{n,2}').$$
By 5.3 above, we have that ${c'}^*(\Sigma U_{n,2}')=0\in H^{2n-3}(T(\bar\sigma_n'))$.
Thus the composite,
\bdm\CD T(\bar\sigma_n')@>{c'}>>
\Sigma T(\sigma_{n,2}')@>{\Sigma U_{n,2}'}>> K_{2n-3},
\endCD \edm
is homotopic to the constant. Then, by taking the suspension, the following is homotopic to
the constant:
\bdm \CD \Sigma T(\bar\sigma_n')@>{\Sigma c'}>> \Sigma^2 T(\sigma_{n,2}')@>\Sigma^2 U_{n,2}'>>\Sigma K_{2n-3}.\endCD \edm

In view of (5.5), this is the composite,
\bdm
\CD
T(\bar\sigma_n)@>c>> \Sigma T(\sigma_{n,2})@>\Sigma U_{n,2}>>\Sigma K_{2n-3}.
\endCD
\edm
\ep

Now we provide:

\vspace{8pt}
\noindent {\it Proof of Theorem B.}\
Let
$\bar\gamma$ denote the universal bundle, whose $(2n-2)$-th Wu class
vanishes, as was introduced in the first paragraph of \S 3 above.
Since the tangent bundle $\bar\tau_n$ of $Z_n$ is stably the pull-back of
$2n\xi_{n-1}$ by the map $\bar p: Z_n\rightarrow P^{n-1}$, any characteristic class
of $\bar\tau_n$ and its normal bundle $\bar\nu_n$ exists in $\bar p^*H^*(P^{n-1})$.
Therefore $v_i(\bar\nu_n)$ vanishes if $n-1\leq i \leq 2n-2$. In particular,
we have $v_{n-1}(\bar\nu_n)=0$ by the condition $\chi(n)=0$.
It follows that there is a symmetric lift
(see \cite[p.103-104]{suth}),
$$a:\bar\nu_n\times\bar\nu_n\rightarrow \bar\gamma.$$

Recall the trivialization (5.4), $\rho':(k-1)\varepsilon\rightarrow \sigma_{n,2}\oplus\nu_{n,2}$,
which is the restriction of (5.3),
$\rho:(k-1)\varepsilon\rightarrow \bar\sigma_n\oplus\bar\nu_n$,
by means of the inclusion $X_{n,2}\hookrightarrow Z_n$.
It is clear that
$a\circ (i_*\times i_*):\nu_{n,2}\times\nu_{n,2}\rightarrow \bar\gamma$ is also symmetric,
where $i_*:\nu_{n,2}\rightarrow \bar\nu_n$
is the bundle map covering the inclusion $X_{n,2}\rightarrow Z_n$.

Let $\bar g:X\rightarrow \Sigma^l T(\bar\sigma_n)$ be the dual of
         $$T(a\circ\bar\Delta): T(2\bar\nu_n)\longrightarrow T(\bar\gamma).$$
Here $\bar\Delta$ denotes the natural bundle map, $2\bar\nu_n\rightarrow \bar\nu\times\bar\nu$,
covering the diagonal $\Delta:Z_n\rightarrow Z_n\times Z_n$.
Also consider the dual $g:X\rightarrow \Sigma^{l+1}T(\sigma_{n,2})$ of $T(2\nu_{n,2})\rightarrow T(\bar\gamma)$,
where $l$ is an appropriate integer.
Note that the map $T(2\nu_{n,2})\rightarrow T(\bar\gamma)$ can be decomposed as
   \bdm \CD T(2\nu_{n,2})@>{T(2i_*)}>>T(2\bar\nu_n)@>T(a\circ\bar\Delta)>> T(\bar\gamma).
   \endCD
   \edm
Since the collapse map $c:T(\bar\sigma_n)\rightarrow \Sigma T(\sigma_{n,2})$
is the dual of $T(2i_*):T(2\nu_{n,2})\rightarrow T(2\bar\nu_n)$, when the dualities are given
by $\rho'$ (see (4.1) above) and by $\rho$ (see (4.2) above), the map
$g:X\rightarrow \Sigma^{l+1}T(\sigma_{n,2})$ can be decomposed as
\bdm
\CD X@>{\bar g}>>\Sigma^l T(\bar\sigma_n)
@>{\Sigma^l c}>>\Sigma^{l+1}T(\sigma_{n,2}).
\endCD
\edm

Note that the map
\bdm
\CD(\Sigma U_{n,2})\circ c : T(\bar\sigma_n) @>>> \Sigma K_{2n-3}\endCD
\edm
is homotopic to the constant by Lemma 5.4.
Now we have that the functional Steenrod square $Sq^{2n-2}_f(\Sigma^{l+1}\iota)$
vanishes since
  $$f=(\Sigma^l((\Sigma U_{n,2})\circ c))\circ \bar g: X\rightarrow \Sigma^{l+1}K_{2n-3}$$
is homotopic to the constant(see (3.3) and also \cite[Corollary 15.10]{steenrod}). Thus we have that
\bdm b_B(\sigma_{n,2},\rho')=Sq^{2n-2}_f(\Sigma^{l+1}\iota)=0\in H^{l+4n-5}(X).\tag{5.6}\edm

By combining (5.6) with Proposition 2.4 and Theorem 3.1,
we have that $\sigma_{n,2}\approx \tau_{n,2}$.
\ \ \ \ \ \ \ \ \ \ \ \ \ \ \ \ \ \ \ \ \ \ \ \ \ \ \ \ \ \ \ \ \ \ \ \ \ \ \ \ \ \ \ \ \
 \ \ \ \ \ \ \ \ \ \ \ \ \ \ \ \ \ \ \ \ \ \ \ \ \ \ \ \ \ \ \ \ \ \ \ \ \ \ \ \ \ \ \ \ \ \ \ \ \ \ \ \ \ \ \ \ \ \ \ $\square$
\vspace{8pt}

It seems worth to note again that $b_B(\sigma_{n,2}, \rho')=0$ holds due to our choice of $\rho'$
(see Remark 4.2 and (5.4) above).

Recall that the tangent bundle $\tau_{n-1}$ of the projective space $P^{n-1}$ is stably
$n\xi_{n-1}$, and therefore the stable inverse of $2n\xi_{n-1}$ is $2\nu_{n-1}$, where
$\nu_{n-1}$ is the normal bundle of $P^{n-1}$. Now we observe the following:

%lemma 5.5
\begin{lemma}
If $n$ is even, we have that $\chi(n)=0$.
\end{lemma}
\bp
Let $v$ denote the total Wu class $v(\nu_{n-1})=c_0+c_1t+\cdots +c_{n-1}t^{n-1}$, where $c_i$ are mod 2 integers, and $c_0=1$.
The odd degree terms vanish for $v^2=v(2\nu_{n-1})$: The coefficient of $t^k$, $0\leq k \leq n-1$, in $v^2$
is given by $\sum_{i+j=k}c_ic_j$. Now let $k$ be odd. Then, if the pair $(i,j)$ is such that $i+j=k$, then $(j,i)$
is also such and $(i,j)\neq (j,i)$. Thus $\sum_{i+j=k}c_ic_j$ vanishes pairwise.
Therefore $\chi(n)$ vanishes, which is the coefficient of $t^{n-1}$ in $v^2$.
\ep

This lemma proves the following, which implies (b), Theorem A above.
\begin{corollary}
If $n$ is even, we have that
$\sigma_{n,2}\approx \tau_{n,2}$.
\end{corollary}

The following is a table for $\chi(n)$ for a few small odd $n$.
\bdm
 \begin{tabular}{|c|c|c|c|c|c|c|c|c|c| }
 \hline
$n$       &  3 & 5 & 7  & 9  & 11 & 13 & 15 & 17 &19 \\
 \hline
$\chi(n)$&  1 & 0 & 0  & 0  & 1  & 0  & 0  & 0  &0 \\
 \hline
\end{tabular}
\tag{5.7}
\edm
Here $\chi(n)$ is calculated as follows. The total Wu class $v(\tau_{n-1})$
of the tangent bundle $\tau_{n-1}$ of $P^{n-1}$ is given by the formula (cf.\ \cite[p.134]{milnor}):
\bdm v(\tau_{n-1})= \sum_{i=0}^{\frac{n-1}{2}}\binom{n-i-1}{i} t^i ,\edm
where $t\in H^1(P^{n-1})$ is the generator of the cohomology algebra. For this calculation, it was convenient for us to use the table on p.~46, \cite{milnor}.
Then the coefficient of $t^{\frac{n-1}{2}}$ in $v(\tau_{n-1})^{-1}$ is $\chi(n)$.

The table implies that $\sigma_{n,2}\approx \tau_{n,2}$ for odd $n$, $3 \leq n \leq 19$,
except possibly
for $n=3, 11$.
In fact, we have $\sigma_{3,2}\approx \tau_{3,2}$ by Lemma 2.2 above and the lemma is consistent
with the above calculation when $n=5$.


\begin{thebibliography}{30}

\bibitem{adams2}
Adams,~J., Vector fields on spheres, Ann.~of Math.~75 (1962), 603-632.

\bibitem{adams3}
Adams,~J., Infinite Loop Spaces, Princeton Univ.~Press, Princeton, N.~J.~1978.

\bibitem{antoniano}
Antoniano,~E., La acci\'on del \'algebra de Steenrod sobre las variedades de Stiefel proyectivas, Bol.~Soc.~Mat.~
Mexicana 22 (1977), 41-47.

\bibitem{antoniano2}
Antoniano,~E., Gitler,~S., Ucci,~J., Zvengrowski,~P.,
On the K-theory and parallelizability of projective Stiefel manifolds, Bol.~Soc.~Mat.~Mexicana 31 (1986), 29-46.

\bibitem{barufatti1}
Barufatti,~N., Obstructions to immersions of projective Stiefel manifolds,
Contemporary Mathematics (Amer.~Math.~Soc.) 161 (1994), 281-287.

\bibitem{barufatti2}
Barufatti,~N., Hacon,~D., K-theory of projective Stiefel manifolds, Trans.~Amer.~Math.~Soc. 352 (2000), 3189-3209.

\bibitem{borel}
Borel,~A., Sur l'homologie et la cohomologie des groupes de Lie compacts connexes, Amer.~J.~Math.~76 (1954), 273-342.

\bibitem{browder}
Browder,~W., The Kervaire invariant of framed manifolds and its generalization,
Ann.~of~Math.~90 (1969), 157-186.

\bibitem{byun1}
Byun,~Y., Tangent fibration of a Poincar\'e complex, J.~London.~Math.~Soc. (2) 59 (1999), 1101-1116.

\bibitem{byun2}
Byun,~Y., Tangent bundle of a manifold and its homotopy type, J.~London.~Math.~Soc. (2) 60 (1999), 301-307.

\bibitem{dupont}
Dupont,~J., On homotopy invariance of the tangent bundle II, Math.~Scand.~26 (1970), 200-220.

\bibitem{Git-Hand}
Gitler,~S., Handel,~D., The projective Stiefel manifolds - I, Topology 7 (1968), 39-46.

\bibitem{james-thomas}
James,~I., Thomas,~E., An approach to the enumeration problem for non-stable vector bundles,
J.~Math.~Mech.~14 (1965), 485-506.

\bibitem{julius-peter}
Korba\v s,~J., Zvengrowski,~P., The vector field problem:
a survey with emphasis on specific manifolds, Exposition Math. 12 (1994), 3-30.

\bibitem{julius-peter3}
Korba\v s,~J., Zvengrowski,~P., On sectioning tangent bundles and other vector bundles,
Proc.~Winter School Geometry and Physics, Srn\'i (Czech Republic) 1994, Rend.~Circ.~Mat.~Palermo (II),
Supplemento 39 (1996), 85-104.

\bibitem{julius-peter2}
Korba\v s,~J., Zvengrowski,~P., The vector field problem for projective Stiefel manifolds,
Bol.~Soc.~Math.~Mexicana (3) Vol.~15 (2009), 219-234.


\bibitem{koschorke}
Koschorke,~U., Vector Fields and Other Vector Bundle Morphisms-A Singularity Approach,
Lecture Notes in Math.~847, Springer-Verlag, Berlin 1981.

\bibitem{lam-1}
Lam,~K.Y., Sectioning vector bundles over real projective space,
Quart.~J.~Math.~Oxford (2) 23 (1972), 97-106.

\bibitem{lam}
Lam,~K.Y., A formula for the tangent bundle of flag manifolds and related manifolds,
Trans.~Amer.~Math.~Soc.~213 (1975), 305-314.

\bibitem{milnor}
Milnor,~J.W., Stasheff,~J.D., Characteristic Classes, Ann.~of Math.~Studies, No.~76.
Princeton Univ.~Press, Princeton, New Jersey 1974.

\bibitem{sankaran-peter} Sankaran,~P., Zvengrowski,~P., Upper bounds for the span of projective Stiefel manifolds, Contemp.~Math. Vol.~407 (2006), 173-183.

\bibitem{steenrod}
Steenrod,~N., Cohomology invariants of mappings, Ann.~of~Math. 50 (1949), 954-989.

\bibitem{stong}
Stong,~R., Semi-characteristics and free group actions, Compositio Math.~29 (1974), 223-248.

%[Sutherland(1976)]
\bibitem{suth}
Sutherland,~W., The Browder-Dupont invariant, Proc.~London Math.~Soc.~(3) 33 (1976), 94-112.

\bibitem{thomas}
Thomas,~E., Vector fields on manifolds, Bull.~Amer.~Math.~Soc.~75 (1969), 643-683.

\bibitem{peter}
Zvengrowski, P., \"Uber die Parallelisierbarkeit von Stiefel Mannigfaltigkeiten, Forschungsinstitut
f\"ur Mathematik, ETH Z\"urich und University of Calgary. April, 1976.




\end{thebibliography}
\end{document}